\documentclass[oneside]{amsart}
\usepackage[utf8]{inputenc}
\usepackage{amsmath}

\title{Integrable sub-Riemannian geodesic flows on the special orthogonal group}
\author[Bravo-Doddoli]{Alejandro Bravo-Doddoli}
\address{Alejandro Bravo-Doddoli;
	{University of Michigan, 530 Church St, Ann Arbor, MI 48109, United States
		\\
		\href{abravodo@umich.edu}{abravodo@umich.edu}}
}

\author[Arathoon]{Philip Arathoon}
\address{Philip Arathoon;
	{University of Michigan, 530 Church St, Ann Arbor, MI 48109, United States
		\\
		\href{philash@umich.edu}{philash@umich.edu}}
}

\author[Bloch]{Anthony M. Bloch}
\address{Anthony M. Bloch;
	{University of Michigan, 530 Church St, Ann Arbor, MI 48109, United States
		\\
		\href{abloch@umich.edu}{abloch@umich.edu}}
}

\date{August 2025}
\usepackage[english]{babel}
\usepackage[utf8]{inputenc}
\usepackage{libertine}
\usepackage{graphicx}
\usepackage{amsthm}
\usepackage{amsmath}
\usepackage{hyperref}
\usepackage{tikz}
\usepackage[all]{xy}
\usepackage[normalem]{ulem} %strikethough with \soutt, just for editing
\usepackage{amsfonts}

\usepackage {amssymb}
\newtheorem{Theorem}{Theorem}
\newtheorem{Prop}{Proposition}[section]

\theoremstyle{remark}
\newtheorem{Remark}{Remark}[section]

\theoremstyle{plain}
\DeclareMathOperator{\Tr}{Tr}
\newcommand{\R}{\mathbf{R}}
\newcommand{\C}{\mathbf{C}}
\newcommand{\D}{\mathcal{D}}
\DeclareMathOperator{\cn}{cn}
\DeclareMathOperator{\sn}{sn}
\DeclareMathOperator{\rank}{\textrm{rank}}

\begin{document}
	
	%\textcolor{blue}{We \sout{analyze} analyse}
	
	\begin{abstract}
		We analyse the geometry of the rubber-rolling distribution on the special orthogonal group and show that almost all the normal geodesics of any right-invariant sub-Riemannian metric defined on this distribution are completely integrable. Our argument is an adaptation of the method used to establish integrability of the Riemannian metric arising from the  \(n\)-dimensional rigid body: namely, by exhibiting a Lax pair and bi-Hamiltonian structure for the reduced equations of motion.

	\end{abstract}
	
	\maketitle
	
	\setcounter{tocdepth}{1}%%%%edited depth
	%	\tableofcontents
	
	\section{Introduction}
	
	The geodesic flow on a Riemannian manifold \(M\) is generated by a Hamiltonian system on the cotangent bundle. A Riemannian metric is called integrable if the corresponding Hamiltonian system is integrable in the sense of Liouville, that is, if the flow admits \(n=\dim M\) independent integrals of motion which mutually commute with respect to the Poisson bracket on \(T^*M\). In principle, one can then find the geodesics explicitly through quadrature. The literature surrounding integrable Riemannian metrics is vast (see, for instance, \cite{double_bracket,thimm,paternain,bolsinov,bolsinov_obstruction}) whereas less has been said about sub-Riemannian metrics, with much of the work focussed on Carnot groups \cite{kruglikov2017integrability,monroy2003integrability,sympletic-reduction}. A sub-Riemannian metric also defines a Hamiltonian system on the cotangent bundle, and so one can equally ask if this flow is integrable. 
	
	A famous example of an integrable Riemannian metric on \(\mathbf{SO}(n)\) is that which arises from the kinetic energy of an \(n\)-dimensional rigid body. We will consider a related class of sub-Riemannian metrics defined on the right-invariant `rubber-rolling' distribution \(\mathcal{D}\) of tangent vectors to \(\mathbf{SO}(n)\) of the form \(\dot{g}=\Omega g\) where
	\begin{equation}\label{rolling_dist}
		\Omega=\left(
		\begin{array}{c|c}
			0 & \begin{matrix}
				\Omega_{12} & \dots & \Omega_{1n}
			\end{matrix}\\\hline
			\begin{matrix}
				-\Omega_{12}\\\vdots\\-\Omega_{1n}
			\end{matrix}
			&  \mbox{\huge $0$}
		\end{array}
		\right)
	\end{equation}
	are skew-symmetric matrices whose only non-zero entries are in the first row and column.

	\begin{Theorem}\label{sr_thm}
		The sub-Riemannian metric with length
		\begin{equation}\label{sr_length}
			\sqrt{A_2\Omega_{12}^2+\dots +A_n\Omega_{1n}^2}
		\end{equation}
		defined on the distribution \(\mathcal{D}\) defines an integrable Hamiltonian system on \(T^*\mathbf{SO}(n)\) for any choice of distinct positive scalars \(A_2,\dots,A_n\).
	\end{Theorem}
	
	%\sout{We show this theorem follows from}{\color{blue}
		
		We prove this theorem by adapting the same methods which led Manakov to establish the complete integrability of the rigid body \cite{manakov}. In particular, we consider the symplectic reduced space \(\mathfrak{so}(n)^*\) and show that the corresponding reduced equations of motion admit a Lax pair formulation and a bi-Hamiltonian structure. This yields a collection of first integrals which can then be shown to define an integrable system thanks to the work of Mishchenko and Fomenko in \cite{misch_fomenko}. The idea of applying Mishchenko and Fomenko's arguments to a sub-Riemannian metric instead of an ordinary Riemannian metric was pointed out in \cite{bloch1994sub} and \cite{singular_manakov}.

			Any sub-Riemannian problem may be reinterpreted as a kinematic non-holonomic optimal control problem. From this point of view, our problem is a special case of a more general class of problems defined on symmetric spaces considered in \cite{bloch1994sub,bloch2015nonholonomic}. There is also a wider literature concerning the related problem of rolling manifolds along each other, see for instance \cite{jurdjevic_conferencepaper,leite}. In particular, we note that Theorem~\ref{sr_thm} was proven in the symmetric case where \(A_2=\dots=A_n\) for \(n=3\) in \cite{Baillieulthesis} and more generally in \cite{jurdjevic_rolling} where the sub-Riemannian problem corresponds to a ball rolling on a plane without spinning.

			\section{A sub-Riemannian metric on the special orthogonal group}\label{sec:pre}
			\subsection{The Hamiltonian setting}
			A sub-Riemannian metric on a manifold \(M\) is completely determined by a fibrewise self-adjoint map \(\beta\colon T^*M\rightarrow TM\) of constant rank known as the cometric. The distribution \(\mathcal{D}\) is the image of \(\beta\) and the inner product between  \(\beta(p_1)\) and \(\beta(p_2)\) is defined as \(\langle p_1,\beta(p_2)\rangle\). The function
			\[
			H(p)=\frac{1}{2}\langle p,\beta(p)\rangle
			\]
			on \(T^*M\) is the Hamiltonian for the sub-Riemannian metric. It is known that the integral curves of the Hamiltonian flow when projected down to \(M\) are, for sufficiently small arcs, the unique length-minimising admissible curves between their endpoints.
			
			We shall be interested in the case when the manifold is a Lie group \(G\). Left and right multiplication of the group on itself lifts to an action on the cotangent bundle, and we will say that the sub-Riemannian metric on \(G\) is left- or right-invariant if the Hamiltonian is invariant with respect to the left or right action. 
			
			The momentum maps
			\begin{equation*}\label{legs}
				\xymatrix{
					& T^*G   \ar[dl]_{\mathbf{J}_L}  \ar[dr] ^{\mathbf{J}_R }&  \\
					%\ar[dr] & \\ 
					{\mathfrak{g}^* _+}   & &  {\mathfrak{g}^* _-} 
				}
			\end{equation*}
			for the left and right actions are given by right and left translation of covectors back to the identity, respectively. For left-invariant Hamiltonians, the momentum \(\mathbf{J}_L\) is constant along the Hamiltonian flow, and \(\mathbf{J}_R\) serves as an orbit-quotient; and vice versa with the words left and right interchanged. The orbit space \(\mathfrak{g}_\pm^*\) is Poisson with the bracket 
			\begin{equation}\label{poisson_bracket}
					\{f,g\}(M)=\pm\langle M,[df|_M,dg_M]\rangle
				\end{equation}
				evaluated at \(M\in\mathfrak{g}^*_\pm\). Here the derivatives \(df\) and \(dg\) evaluated at \(M\) are canonically identified with elements of \(\mathfrak{g}\) and \([~,~]\) denotes the Lie bracket. 
			
			For a more detailed background in sub-Riemannian geometry, we recommend the references \cite{tour}, \cite{comprehensiveSR}, and \cite{bloch2015nonholonomic}, and \cite{ortega_ratiu} for those aspects concerning Poisson geometry.

			\subsection{The rigid body}
			The kinetic energy for an \(n\)-dimensional rigid body defines a left-invariant Riemannian metric on \(\mathbf{SO}(n)\). This descends through reduction to \(\mathfrak{so}(n)^*_-\) with Hamiltonian
			\[
			H(M)=\frac{1}{2}\langle M,\Omega\rangle,\quad\text{for}~\Omega=\mathbb{J}^{-1}(M)
			\]
			where \(\mathbb{J}\colon\mathfrak{so}(n)\rightarrow\mathfrak{so}(n)^*\) is the inertia operator \cite{tudor-free-rigid-body}. If we identify \(\mathfrak{so}(n)\) with its dual via the trace pairing 
				\[\langle M,\Omega\rangle=\frac{1}{2}\Tr (M^T\Omega)\]
			then
			\[
			\mathbb{J}(\Omega)=J\Omega+\Omega J
			\]
			where \(J\) is the mass matrix, which we may suppose is diagonal with strictly positive entries \(J_1,\dots, J_n\). The angular momentum is related to the angular velocity via
			\[
			M_{kl}=I_{kl}\Omega_{kl},\quad\text{for}~I_{kl}=J_k+J_l.
			\]
			The \(I_{kl}\) are the generalised moments of inertia for the \(n\)-dimensional body, and the Hamiltonian is in coordinates given by
			\begin{equation}\label{rigid_body_ham}
				H(M)=\frac{1}{2}\sum_{k<l}\frac{M_{kl}^2}{I_{kl}}.
			\end{equation}
			
			The equations of motion on \(\mathfrak{so}(n)^*_-\) are known as the Euler equations, and given by
			\[
			\frac{d}{dt}M=[M,\Omega].
			\]
			The trivial constants of motion are the Hamiltonian and the Casimirs \(\Tr(M^k)\). Additional constants of motion can be seen by rewriting the Euler equations in terms of a Lax pair
			\[
			\frac{d}{dt}(M+\lambda J^2)=[M+\lambda J^2,\Omega+\lambda J]
			\]
			where \(\lambda\) is an arbitrary parameter. This equation, found originally by Manakov \cite{manakov}, implies that solution curves remain inside the same conjugacy class, and so \(\Tr(M+\lambda J^2)^k\) is constant for any power \(k\). Since \(\lambda\) was arbitrary it follows that each of the coefficients of \(\lambda\) in the expansion
			\begin{equation}\label{manakov_integrals}
				h^\lambda_k=\Tr(M+\lambda J^2)^k=\sum_{l=0}^{k}\lambda^{k-l}h_{k,l}(M),
			\end{equation}
			is itself a constant of motion. These define the {Manakov integrals} \(h_{k,l}\).

			\subsection{The sub-Riemannian metric}
			Let \(\mathfrak{p}\) denote the subspace of matrices of the form \eqref{rolling_dist} and equip this with the metric given in \eqref{sr_length}. We note that this is not a Lie subalgebra of \(\mathfrak{so}(n)\), but that it generates the entire algebra in one step, that is \[\mathfrak{p}+[\mathfrak{p},\mathfrak{p}]=\mathfrak{g}.\] 
			Left/right multiplication pushes \(\mathfrak{p}\) forward to define a non-integrable left/right-invariant distribution on \(\mathbf{SO}(n)\) equipped with a sub-Riemannian metric.  The resulting length functional on the left-invariant distribution corresponds to the \(n\)-dimensional Suslov problem which concerns the motion of a rigid body with constrained angular momentum.  However, the dynamics of this system are non-holonomic, and as explained in \cite{zb00} do not correspond to a system of Lie-Poisson equations on \(\mathfrak{so}(n)_-^*\) but instead to a system of so-called Euler-Poincar\'{e}-Suslov equations. Therefore, we focus our attention instead on the right-invariant distribution whose length functional may be interpreted as describing the length traced by the path of an \((n-1)\)-dimensional ball which rolls without sliding or spinning on \(\R^{n-1}\) with metric 
			\begin{equation}\label{hyperplane_metric}
				ds^2=A_2dy_2^2+\dots+A_ndy_n^2.
			\end{equation}
			We denote this right-invariant distribution by \(\mathcal{D}\) and refer to it as the rubber-rolling distribution, a term we have taken from \cite{ekm} and \cite{rubberrolling}, so called since the ball is not allowed to spin about its point of contact. We will later solve for the rolling motion in the case \(n=3\).
			
			The right-invariant sub-Riemannian Hamiltonian on \(T^*\mathbf{SO}(n)\) descends to give a Lie-Poisson system on \(\mathfrak{so}(n)^*_+\) whose Hamiltonian is given in coordinates by
			\begin{equation}\label{sr_ham}
				H_\text{sR}(M)=\frac{1}{2}\sum_{k=2}^n\frac{M_{1k}^2}{A_k},
			\end{equation}
			where, as we did earlier, we have identified the Lie algebra with its dual by using the trace form. 
			
			\begin{Remark}\label{limit_remark}
				Notice that the sub-Riemannian Hamiltonian in \eqref{sr_ham} is the limit of the rigid-body Hamiltonian in \eqref{rigid_body_ham} when \(I_{kl}\rightarrow\infty\) for \(k\ne 1\) and \(A_k=I_{1k}\). Indeed, any sub-Riemannian metric on a distribution \(\D\) can be interpreted as a limit of Riemannian penalty metrics in which velocity vectors not belonging to \(\D\) become increasingly expensive. We note, however, that this limit is not physical since \(2J_1=I_{12}+I_{13}-I_{23}\) necessarily limits to \(-\infty\), which contradicts the positive-definiteness of the mass matrix for realistic bodies.  On the other hand, the metric on the left-invariant distribution does correspond to the optimal control of an underactuated $n$-dimensional rigid body, with control torques (actuators) about $(n-1)$ of the $n$ principal axes.  For the underactuated rigid body problem, see e.g. \cite{bloch1994sub} and \cite{bloch2015nonholonomic} and references therein. 
			\end{Remark}	
			\subsection{Rolling a ball across a table}\label{sec:rolling}
			
			Place a rubber ball on a table. The configuration of the ball is determined by its point of contact \((y,z)\in\R^2\) and a matrix \(g\in\mathbf{SO}(3)\) whose columns are the basis vectors of an orthonormal frame fixed inside the ball. How can we roll the ball from one orientation \(g_0\) to another \(g_1\) in such a way that the curve of contact which the ball traces on the table has the least length with respect to the metric
			\begin{equation*}\label{table_metric}
				ds^2=A_2dy^2+A_3dz^2?
			\end{equation*}
			This question is equivalent to our sub-Riemannian problem on \(\mathbf{SO}(n)\) for \(n=3\), as we will now show.
			
			The no-spinning constraint implies that admissible tangent vectors to \(\mathbf{SO}(3)\) must be of the form 
			
				\begin{equation}\label{distribution}
					\frac{dg}{dt}=\begin{pmatrix}
						0&-\Omega_{z}&\Omega_{y}\\
						\Omega_{z}& 0 & 0\\
						-\Omega_{y} & 0 & 0
					\end{pmatrix}g=\widehat{\Omega}g.
				\end{equation}
			
			This notation is chosen to coincide with the `hat-map' whereby \(\widehat{\Omega}v=\Omega\times v\) for \(\Omega=(0,\Omega_y,\Omega_z)^T\). For such a motion, the path traced on the plane satisfies
				\begin{equation}\label{path}
					\frac{d}{dt}\begin{pmatrix}
						y\\z
					\end{pmatrix}=\begin{pmatrix}
						\Omega_{z}\\ -\Omega_{y}
					\end{pmatrix},
				\end{equation}
				where we are supposing for simplicity that the ball has a unit radius. For an admissible curve \(g(t)\) the length of the path traced by the ball is therefore
			\begin{equation}\label{SR_metric}
					\int \sqrt{A_2|\Omega_{z}|^2+A_3|\Omega_{y}|^2}~dt.
				\end{equation}
			Eq.~\eqref{distribution} defines a right-invariant distribution on \(\mathbf{SO}(3)\) and Eq.~\eqref{SR_metric} equips this with a sub-Riemannian metric.
			
			\begin{figure}
				\begin{center}
					\begin{tikzpicture}
						\draw (-3,4) node[inner sep=0] {\includegraphics[scale=0.6]{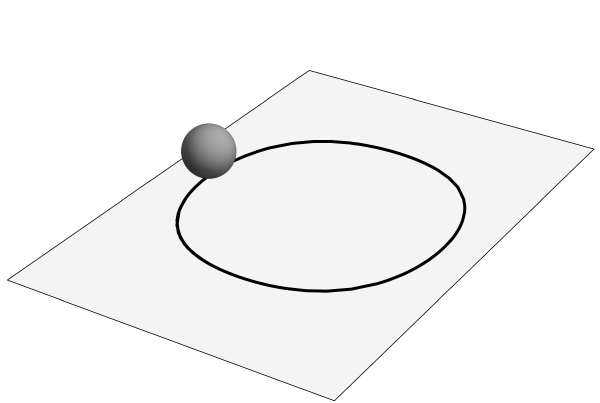}};
						\draw (-3,-1) node[inner sep=0] {\includegraphics[scale=0.6]{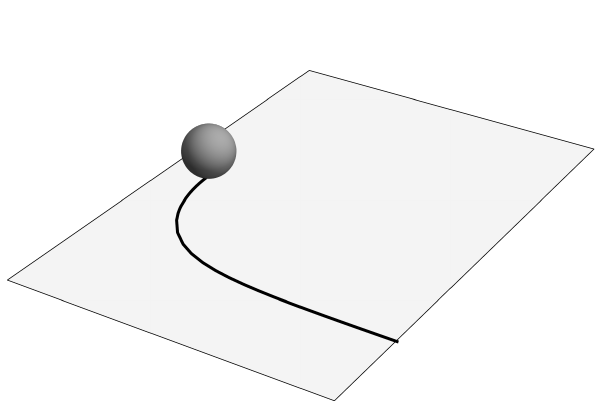}};
						\draw (-3,-6) node[inner sep=0] {\includegraphics[scale=0.6]{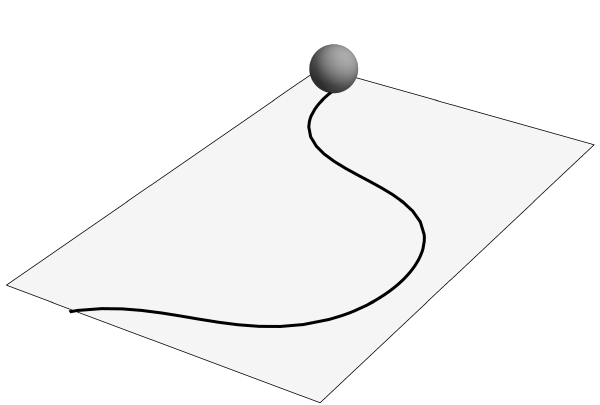}};
						
						\draw (3,4-.5) node[inner sep=0] {\includegraphics[scale=0.6]{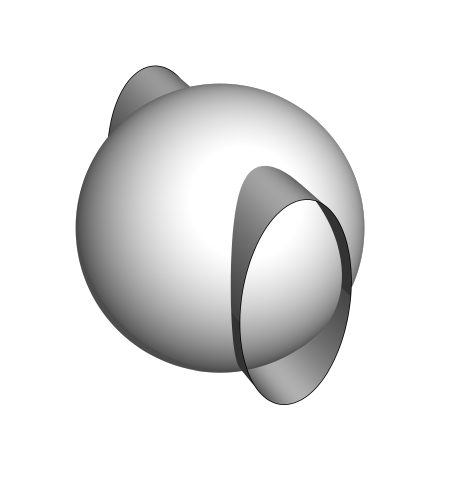}};
						\draw (3,0-.5-1) node[inner sep=0] {\includegraphics[scale=0.6]{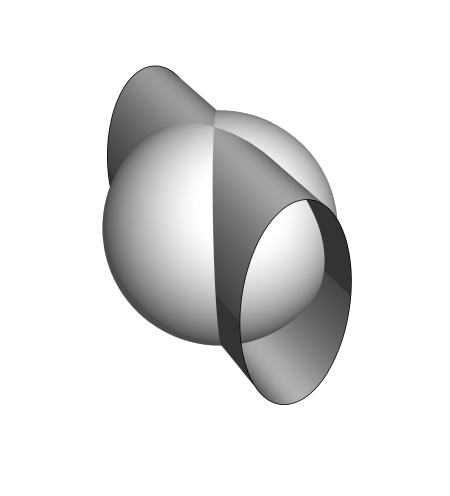}};
						\draw (3,-4-.5-2) node[inner sep=0] {\includegraphics[scale=0.6]{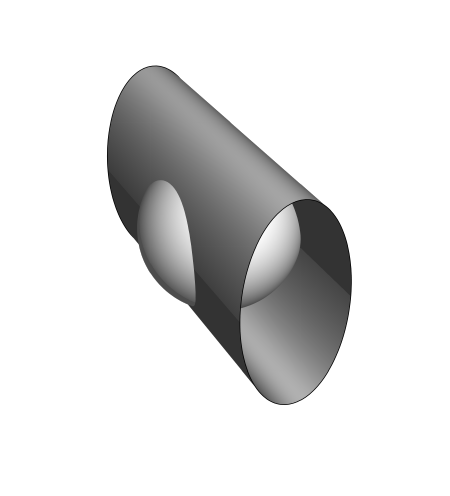}};
						
						\draw (-6,6) node {$0\le k<1$};
						\draw (-6,1) node {$k=1$};
						\draw (-6,-4) node {$k>1$};
					\end{tikzpicture}
				\end{center}
				\caption{\label{rolling_balls}On the right are intersections of the elliptic-cylinder $H_\text{sR}=1/2$ with coadjoint orbits in $\mathfrak{so}(3)^*$. The corresponding paths which the ball traces on the table are shown on the left. }
			\end{figure}
			
			Identify tangent vectors \(\dot{g}\) to \(\mathbf{SO}(3)\) with vectors \(\Omega\) in \(\R^3\) via right-trivialisation \(\widehat{\Omega}=\dot{g}g^{-1}\) and the hat-map. Likewise, we may also identify covectors with \(M\) in \(\R^3\) via the standard pairing \(\langle M,\Omega\rangle=M^T\Omega\). The Hamiltonian on \(T^*\mathbf{SO}(3)\) for the sub-Riemannian metric is then
				\begin{equation*}
					H_\text{sR}(M)=\frac{1}{2}\left(\frac{M_{z}^2}{A_2}+\frac{M_{y}^2}{A_3}\right)
				\end{equation*}
			for \(M=(M_{x},M_{y},M_{z})^T\). After reduction this descends to a Lie-Poisson system in \(\mathfrak{so}(3)^*_+\cong\R^3\) with equation of motion
			\begin{equation}\label{Lie-Poisson}
				\frac{d}{dt}M=\nabla H_\text{sR}\times M.
			\end{equation}
			The Hamiltonian \(H_\text{sR}\) and Casimir function \(|M|^2\) are both constants of motion. The intersections of their levels sets therefore give the solution curves: elliptic cylinders for \(H_\text{sR}\) constant, and spheres of constant \(C=|M|^2\), as shown on the right-hand side in Figure~\ref{rolling_balls}. 
			
			For the initial conditions \(M_{x}(0)=\sqrt{C-A_2}\), \(M_{y}(0)=0\) , and \(M_{z}(0)=\sqrt{A_2}\) we can solve \eqref{Lie-Poisson} exactly in terms of Jacobi elliptic functions, finding
			\begin{equation}\label{Jacobi}	M_{y}(t)=\sqrt{A_3}\sn(u,k)\quad\text{and}\quad M_{z}(t)=\sqrt{A_2}\cn(u,k),
				\end{equation}
			where
			\begin{equation*}
				u=t\sqrt{\frac{C-A_2}{A_2A_3}}\quad\text{and}\quad k=\sqrt{\frac{A_3-A_2}{C-A_2}}.
			\end{equation*}
			The tangent vector \(\Omega\) is related to the covector \(M\) via \(M_{y}=A_3\Omega_{y}\) and \(M_{z}=A_2\Omega_{z}\). We can use this relation together with \eqref{path} to find the path that the ball traces by integrating \eqref{Jacobi} (with the aid of \cite{Jacobi}) to find that
			\begin{equation}\label{curves}
				\cos\left(z\sqrt{\textstyle\frac{A_3-A_2}{A_3}}\right)=
				\begin{cases}
					\sqrt{1-k^2}\cosh\left(\alpha-y\sqrt{\frac{A_3-A_2}{A_2}}\right),&0\le k<1\\[20pt]
					\exp\left(-y\sqrt{\frac{A_3-A_2}{A_2}}\right), & k=1\\[20pt]
					\sqrt{k^2-1}\sinh\left(\alpha-y\sqrt{\frac{A_3-A_2}{A_2}}\right),&k>1
				\end{cases}
			\end{equation}
			where we are assuming \(A_3>A_2\) and \(\alpha\) is such that \(y(0)=z(0)=0\). These curves are shown on the left-hand side in Figure~\ref{rolling_balls}.

			\begin{Remark}
				The rolling problem we have discussed is related to an interview question asked in the 1950s at the University of Oxford \cite{hammer}, first solved in \cite{arthur} (see also \cite{plateball} and \cite{jurjevic}). Given two points \(p_0\) and \(p_1\) on the table, and two configurations \(g_0\) and \(g_1\) of the ball, what is the shortest length path which rolls \((p_0,g_1)\) to \((p_1,g_1)\)? Our problem differs from this in two respects. Firstly, we are not concerned with connecting two points on the table, but rather with the two configurations of the ball; and secondly, the metric on the table is not the standard Euclidean metric. It is curious to compare the curves in Eq.~\eqref{curves} with those for the Euclidean case \(A_2=A_3\) where the ball is found to travel across the table in straight lines and circular arcs.
			\end{Remark}
			
			We conclude by noting that we have only solved for the path which the ball traces on the table. In order to find the explicit sub-Riemannian geodesics on the group \(\mathbf{SO}(3)\) Eq.~\eqref{distribution} must be integrated. These geodesics have been found explicitly in \cite{sachkov}.

			\section{Integrability}\label{sec:sR-strucutre}
			\subsection{Manakov integrals in the limit}
			Following on from Remark~\ref{limit_remark} we may interpret the sub-Riemannian metric as a limit of Riemannian metrics for the rigid body. Since the rigid body is known to be integrable for distinct principal moments of inertia, does it follow that its limit is also integrable? Our initial heuristic for finding first integrals was to inspect the Manakov integrals in such a limit.  However, the integrals are found to diverge or become dependent. To overcome this, we now show how to define a collection of integrals which are functionally equivalent to the Manakov integrals, but which survive in the limit. We will later show that these are in fact first integrals of the sub-Riemannian flow.
			
			To define this limit more precisely, consider a rigid-body parametrised by \(s\) whose mass matrix \(J^s\) is diagonal with entries
			\[
			J_1^s=-s,\text{ and }J_k^s=A_k+s\text{ for }k\ne 1.
			\]
			In this way \(I_{1k}=A_k\) for all \(s\), and the \(I_{kl}\rightarrow\infty\) as \(s\rightarrow\infty\) for \(k\ne 1\). The Hamiltonian therefore limits to the sub-Riemannian Hamiltonian \(H_\text{sR}\). Strictly speaking, this metric does not represent a physical body since \(J^s_1\) is negative, however this does not affect the definition of the Manakov integrals. We may rewrite
			\[
			h_k^\lambda=\Tr(M+\lambda \left[J^s\right]^2)^k
			\]
			as
			\[
			\Tr(M+\lambda\left[s^2\text{Id}+2sA+A^2\right])^k,
			\]
			where we have introduced the diagonal matrix \(A=\text{diag}(0,A_2,\dots, A_n)\). Replace the indeterminate variable \(\lambda\) with \(\sigma/(2s)\) so that this becomes
			\[
			h_k^\sigma=\Tr\left(\left[M+\sigma\left(A+\frac{A^2}{2s}\right)\right]+\frac{s\sigma}{2}\text{Id}\right)^k.
			\]
			Notice that the expression above may be written as
			\[
			\Tr\left[M+\sigma\left(A+\frac{A^2}{2s}\right)\right]^k+\text{ a linear combination of }h^\sigma_{k-1},\dots, h^\sigma_2, h^\sigma_1, h^\sigma_0.
			\]
			Therefore, the Manakov integrals \(h_k^\sigma\) are functionally equivalent to the collection
			\[
			\Tr\left[M+\sigma\left(A+\frac{A^2}{2s}\right)\right]^k
			\]
			for \(k=1,2,\dots, n\). These integrals have a well-defined limit as \(s\rightarrow\infty\) given by
			\[
			g^\sigma_k=\Tr(M+\sigma A)^k
			\]
			for arbitrary \(\sigma\). In the next section, we will verify that these are indeed first integrals of the sub-Riemannian geodesic flow.

			\subsection{A Lax pair and bi-Hamiltonian structure}\label{int_section}
			
			The Poisson structure in \eqref{poisson_bracket} leads to the equation of motion
			\begin{equation}\label{Lie_Poisson_big}
					\frac{d}{dt}M=[dH_\text{sR}|_M,M],
				\end{equation}
			in \(\mathfrak{so}_+(n)^*\) where the derivative \(dH_\text{sR}\) evaluated at \(M\) is understood to be an element of the Lie algebra. It will be helpful to decompose 
			\begin{equation}\label{decomp}
				M=\begin{pmatrix}
					0 & -\mu^T\\
					\mu&\widehat{m}
				\end{pmatrix}
			\end{equation}
			into the \((n-1)\)-dimensional column vector \(\mu\) and the skew-symmetric matrix \(\widehat{m}\) belonging to \(\mathfrak{so}(n-1)\). In a slight abuse of notation, by letting \(\hat{A}^{-1}\) denote the matrix \(\text{diag}(A_2^{-1},\dots, A_n^{-1})\) the Hamiltonian is conveniently written
			\[
			H_\text{sR}(M)=\frac{1}{2}\Tr\left(\mu^T\hat{A}^{-1}\mu\right)
			\]
			from which we find
			\begin{equation}\label{sr_ham_derivative}
					dH_\text{sR}|_M=\begin{pmatrix}
						0 & -\mu^T\hat{A}^{-1}\\
						\hat{A}^{-1}\mu& 0
					\end{pmatrix}.
				\end{equation}
			Together with \eqref{Lie_Poisson_big} this allows us to write the equations of motion as
			\begin{equation}\label{decomposed_eqns}
				\begin{split}
					\frac{d}{dt}\mu&=-\widehat{m}\hat{A}^{-1}\mu,\\
					\frac{d}{dt}\widehat{m}&=[\mu\mu^T,\hat{A}^{-1}].
				\end{split}
			\end{equation}

				\begin{Prop}
					The equations of motion admit a Lax pair formulation,
					\begin{equation}\label{sr_lax_pair}
						\frac{d}{dt}(M+\sigma A)=[dH_\textup{sR}|_M+\sigma B,M+\sigma A],
					\end{equation}
					for \(\sigma\) an arbitrary parameter, \(A=\textup{diag}(0,A_2,\dots,A_n)\), and \(B=\textup{diag}(0,1,\dots,1)\).
				\end{Prop}
				\begin{proof}
					This follows from the important identity
					\begin{equation}\label{important}
						[dH_\textup{sR}|_M,A]=[M,B]
					\end{equation}
					which is verified by direct calculation.
				\end{proof}

			It follows that \(M+\sigma A\) evolves within the same conjugacy class. The traces of \({(M+\sigma A)^k}\) are therefore integrals of motion. Since \(\sigma\) is arbitrary, the coefficients of \(\sigma\) in the expansion
			\begin{equation}\label{sr_manakov_integrals}
				g^\sigma_k=\Tr(M+\sigma A)^k=\sum_{l=0}^{k}\sigma^{k-l}g_{k,l}(M)
			\end{equation}
			are first integrals, which we name the sub-Riemannian Manakov integrals, and denote by \(g_{k,l}\).

			\begin{Remark}
				It is interesting to note that the sub-Riemannian Manakov integrals are the same as the ordinary Manakov integrals in \eqref{manakov_integrals} for a `flat' rigid body whose diagonal mass matrix satisfies \(J^2=A\). The body is flat because it has no extent in the first principal direction since \(J_1=0\). The equations of motion for the flat body are, however, not the same as those for the sub-Riemannian problem. This can be seen by noting that the inertia operator \(\mathbb{J}\) for a flat rigid body (with only \(J_1\) equal to zero) is invertible, yet the cometric \(M\mapsto dH_\text{sR}\) for the sub-Riemannian metric is not invertible.
			\end{Remark}
			We now turn our attention to establishing a bi-Hamiltonian structure for the equations of motion. This can be used to give an elegant proof of the involutivity of the sub-Riemannian integrals by closely adapting the work of \cite{bi-hamiltonian}. We will, however, prove this using an alternative method, owing to Mishchenko and Fomenko, in the next section. Nonetheless, we include the work below since it constitutes an interesting example of a sub-Riemannian flow admitting a bi-Hamiltonian structure.
				
				We begin by observing that for any matrix \(S\) there is an alternative Lie bracket on the space of matrices given by 
				\[
				[X,Y]_S=XSY-YSX.
				\]
				If \(S\) is symmetric then this restricts to a Lie bracket on \(\mathfrak{so}(n)\) which we can use to modify \eqref{poisson_bracket} to produce an alternative Poisson structure on \(\mathfrak{so}(n)^*\) given by
				\[
				\{f,g\}_S=\langle M,[df|_M,dg|_M]_S\rangle.
				\]
			
			For a given Hamiltonian \(H\) suppose that there is another function \(F\) for which
			\[
			\{\varphi,H\}=\{\varphi,F\}_S
			\]
			holds for all \(\varphi\). The Hamiltonian system generated by \(H\) is therefore equivalent to the system generated by \(F\) with respect to the altered Poisson structure. If \(S\) is not a scalar multiple of the identity, then we say that this system is bi-Hamiltonian.
			
			\begin{Prop}
				The sub-Riemannian Hamiltonian system generated by \(H_\text{sR}\) on \(\mathfrak{so}_+(n)^*\) is bi-Hamiltonian. Specifically, the system is equivalent to the Hamiltonian system generated by
				\begin{equation}\label{big_F}
					F(M)=\Tr\left(\mu^T\hat{A}^{-1}\mu\right)+\frac{1}{2}\Tr\left(\widehat{m}^T\hat{A}^{-1}\widehat{m}\hat{A}^{-1}\right)
				\end{equation}
				with respect to the altered Poisson structure \(\{~,~\}_A\) where \(A=\textup{diag}(0,A_2,\dots,A_n)\).
			\end{Prop}
			\begin{proof}
				From the definitions of the Poisson brackets, we must show that
				
                    \begin{equation}\label{check_this}
						\langle M,[X,dH_\text{sR}|_M]\rangle=\langle M,[X,dF|_M]_A\rangle
					\end{equation}
				for all \(X\) in \(\mathfrak{so}(n)\). As we are identifying \(\mathfrak{so}(n)\) with its dual via the trace form, one can show that 
					\[\langle M,[X,dH_\text{sR}|_M]\rangle=\langle X,[dH_\text{sR}|_M,M],\rangle\text{ and }
					\langle M,[X,dF|_M]_A\rangle=\langle X,[A,M]_{dF|_M}\rangle.
					\]
				Therefore, it suffices to show 
				\begin{equation}\label{checkthis2}
						\langle X,[dH_\text{sR}|_M,M]\rangle=\langle X,[A,M]_{dF|_M}\rangle.
					\end{equation}
				The derivative of \(F\) is found to be 
					\begin{equation}
						dF|_M=\begin{pmatrix}
							0 & -\mu^T\hat{A}^{-2}\\
							\hat{A}^{-2}\mu & \hat{A}^{-1}\widehat{m}\hat{A}^{-1}
						\end{pmatrix}
					\end{equation}
				which when combined with \eqref{sr_ham_derivative} allows us to compute and verify the important identity \([dH_\text{sR}|_M,M]=[A,M]_{dF|_M}\) from which \eqref{checkthis2} follows.
				
			\end{proof}

			\begin{Remark}
				The bi-Hamiltonian structure we have constructed is a variation of that for the rigid body given in \cite{bi-hamiltonian} (see also the construction in \cite{bloch2005class}) who are careful to note that if any of the \(J_i\) are zero then this ``...can also be employed for flat bodies, but the construction of the Hamiltonian functions is more complicated''. Indeed, we observe that the altered Hamiltonian \(F\) in \eqref{big_F} differs from the form given in \cite[Eq.~3.8]{bi-hamiltonian}. Nevertheless, as confirmed in their concluding remarks, the proof that the integrals are in involution using a recursion relation may still be applied. Indeed, one can reproduce their arguments exactly to show that the sub-Riemannian Manakov integrals commute.
			\end{Remark}
            
		\subsection{The four-dimensional case}
		The case for \(n=4\) is the smallest case to feature a non-trivial sub-Riemannian Manakov integral. The Hamiltonian from Eq~\eqref{sr_ham} is
		\[
		H_\textrm{sR}(M)=\frac{1}{2}\left(\frac{M_{12}^2}{A_2}+\frac{M_{13}^2}{A_3}+\frac{M_{14}^2}{A_4}\right),
		\]
		for \(M\) in \(\mathfrak{so}(4)\). The 6-dimensional Lie algebra \(\mathfrak{so}(4)\) possesses two Casimirs,
		\begin{align*}
			C_1&=M_{12}^2+M_{13}^2+M_{14}^2+M_{34}^2+M_{24}^2+M_{23}^2,\\
			C_2&=M_{12}M_{34}-M_{13}M_{24}+M_{14}M_{23}.
		\end{align*}
		It is convenient to use the decomposition given in Eq.~\eqref{decomp} together with the hat-map to identity \(\widehat{m}\) in \(\mathfrak{so}(3)\) with a vector \(m\) in \(\R^3\). In this way
		\[
		\mu=\begin{pmatrix}
			-M_{12}\\-M_{13}\\-M_{14}
		\end{pmatrix},~
		m=\begin{pmatrix}
			-M_{34}\\M_{24}\\-M_{23}
		\end{pmatrix},
		\]
	and the Casimirs are compactly written
	\begin{align*}
		C_1&=|\mu|^2+|m|^2,\\
		C_2&=\mu^Tm.
	\end{align*}
 
By considering the map \((\mu,m)\mapsto(\mu+m,\mu-m)\), we see that the set of constant \(C_1\) and \(C_2\) is diffeomorphic to the product of two spheres with respective radii \(C_1+2C_2\) and \(C_1-2C_2\). Incidentally, this establishes that the Casimirs satisfy  
\[
C_1\ge 0,\text{ and }C_1^2\ge 4C_2^2.
\]
The regular coadjoint orbits are thus 4-dimensional. Integrability, therefore, requires one further integral in addition to the Hamiltonian. This is given by the sub-Riemannian integral \(g_{3,2}\), which for convenience we scale by \(-\frac{1}{3}\) and denote by
		\[
		F(M)=A_2M_{12}^2+A_3M_{13}^2+A_4M_{14}^2+(A_3+A_4)M_{34}^2+(A_2+A_4)M_{24}^2+(M_2+M_3)M_{23}^2.
		\]
Hence, by Liouville's theorem, \(\mathfrak{so}(4)\) is foliated almost everywhere into tori equal to the connected components of the regular fibres of \((H_\textrm{sR},F,C_1,C_2)\).

Finally, we consider the analogous `rolling motion' as in Section~\ref{sec:rolling} for the case of a 4-dimensional ball rolling without spinning along a hyperplane. We proceed in greater generality, considering a unit ball in \(\R^n\) rolling without slipping along the hyperplane \(y_1=-1\).

Fix an orthonormal frame inside the centre of the ball which is coincident with the standard frame in space at time \(t=0\). Let \(c(t)\) be the vector for the point of contact between the ball and the hyperplane viewed from within the ball's frame. For \(t=0\) this is \(c_0=(-1,0,\dots,0)^T\), and for later \(t\) is given by \(c(t)=g(t)^{-1}c_0\), where \(g(t)\) is the matrix in \(\mathbf{SO}(n)\) whose columns form the ball's frame in space at time \(t\). The velocity of the point of contact viewed from within the ball's frame is 
\[
\dot{c}=-g^{-1}\dot{g}g^{-1}c_0.
\]
To convert vectors from the ball's frame to the space frame, we need only multiply on the left by \(g\). Therefore, the velocity of the contact point in space is \(-\dot{g}g^{-1}c_0\), and thus, the no-slip condition implies that the point of contact satisfies
\begin{equation}\label{y_eqn}
	\frac{d}{dt}y_j=\Omega_{1j},\quad\text{for }2\le j\le n.
\end{equation}
Here we have introduced the right-invariant velocity \(\Omega=\dot{g}g^{-1}\) in \(\mathfrak{so}(n)\). The no-spin condition implies that \(\Omega\) must belong to the rubber-rolling distribution given in Eq.~\eqref{rolling_dist}, and the metric given in Eq.~\ref{hyperplane_metric} defines the sub-Riemannian metric on \(\mathcal{D}\) as in Eq.~\eqref{sr_length}.

\begin{figure}
			\begin{center}
				\begin{tikzpicture}
					\draw (0,0) node[inner sep=0] {\includegraphics[scale=1]{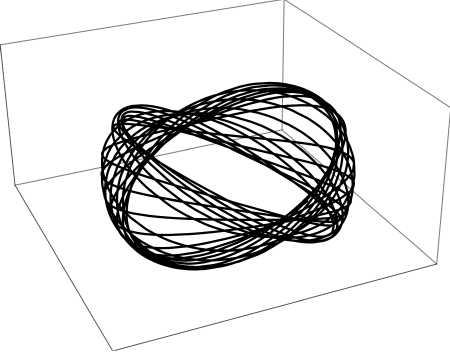}};
					\draw (-4,1) node {$y_4$};
					\draw (-3,-2) node {$y_2$};
					\draw (1.5,-2.5) node {$y_3$};
				\end{tikzpicture}
			\end{center}
			\caption{\label{bounded_curve}A bounded rolling curve.}
		\end{figure}

\begin{figure}
			\begin{center}
				\begin{tikzpicture}
					\draw (0,0) node[inner sep=0] {\includegraphics[scale=1]{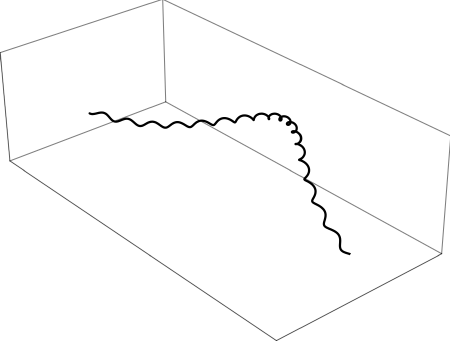}};
					\draw (-4,1) node {$y_4$};
					\draw (-2,-1.5) node {$y_2$};
					\draw (2.5,-2.5) node {$y_3$};
				\end{tikzpicture}
			\end{center}
			\caption{\label{unbound_curve}An unbounded rolling curve.}
		\end{figure}

The equations of motion found in Eq.~\eqref{decomposed_eqns} define the following system in \(\mathfrak{so}(4)\)
		\begin{equation*}
			\arraycolsep=15pt\def\arraystretch{3}
			\begin{array}{ll}
				\displaystyle\frac{dM_{12}}{dt}=-\left(\frac{M_{13}M_{23}}{A_3}+\frac{M_{14}M_{24}}{A_4}\right),&\displaystyle\frac{dM_{23}}{dt}=\left(\frac{A_2-A_3}{A_2A_3}\right)M_{12}M_{13},\\
				\displaystyle\frac{dM_{13}}{dt}=\left(\frac{M_{12}M_{23}}{A_2}-\frac{M_{14}M_{34}}{A_4}\right),&\displaystyle\frac{dM_{24}}{dt}=\left(\frac{A_2-A_4}{A_2A_4}\right)M_{12}M_{14},\\
				\displaystyle\frac{dM_{14}}{dt}=\left(\frac{M_{12}M_{24}}{A_2}+\frac{M_{13}M_{34}}{A_3}\right),&\displaystyle\frac{dM_{34}}{dt}=\left(\frac{A_3-A_4}{A_3A_4}\right)M_{13}M_{14}.
			\end{array}
		\end{equation*}
		The cometric relates the angular velocity with the angular momentum via \(M_{1j}=A_j\Omega_{1j}\) for \(2\le j\le n\). Therefore, we may find the motion for the point of contact of the ball by first integrating the system above, and then integrating Eq.~\eqref{y_eqn}. Examples of these curves for \(n=4\) and for \((A_2,A_3,A_4)=(2,3,4)\) are shown in Figures~\ref{bounded_curve} and \ref{unbound_curve}.
	
		\subsection{Proof of Theorem \ref{sr_thm}}
		The integrability of the sub-Riemannian geodesic system follows directly from an application of  Mishchenko and Fomenko's generalisation of the Euler equations for the rigid body to arbitrary semisimple Lie algebras. For convenience we reproduce their result below \cite[Theorem~5.2]{misch_fomenko} (see also the statement in \cite[Theorem~1]{bolsinov_sectional}).
			\begin{Theorem}\label{mf_thm}
				Let \(\mathfrak{g}_n\subseteq\mathfrak{g}_u\) be the normal subalgebra of the semisimple complex algebra \(\mathfrak{g}\) with compact real form \(\mathfrak{g}_u\) and consider the Euler equations
				\[
				\frac{d}{dt}M=[M,\varphi(M)]
				\]
				defined on the subalgebra \(\mathfrak{g}_n\) where \(\varphi\colon\mathfrak{g}_n\rightarrow\mathfrak{g}_n\) is a self-adjoint linear map satisfying
				\begin{equation}\label{sectional}
					[\varphi(M),A]=[M,B]
				\end{equation}
				for some commuting \(A,B\in\mathfrak{g}\) with \(A\) regular. Then for any Casimir \(f\) of \(\mathfrak{g}_n\) the functions of the form \(f(M+\sigma A)\) are constants of motion for the Euler equations and are in pairwise involution. Furthermore, among this set of functions, one can choose functionally independent integrals equal in number to half the dimension of a regular adjoint orbit of \(\mathfrak{g}_n\).
			\end{Theorem}
            
			To apply this theorem we set \(\mathfrak{g}=\mathfrak{gl}_n\C\) with compact form \(\mathfrak{g}_u=\mathfrak{u}(n)\) and normal form \(\mathfrak{g}_n=\mathfrak{so}(n)\). The sectional operator \(\varphi\) is identified with the map \(M\mapsto dH_\text{sR}|_M\) and the identity \eqref{sectional} is satisfied as in \eqref{important} for \(A=\text{diag}(0,A_2,\dots, A_n)\) and \(B=\text{diag}(0,1,\dots,1)\). Casimirs are invariants of the adjoint representation, which for \(\mathfrak{so}(n)\) are the trace-powers \(f(M)=\Tr M^k\). Therefore, it follows from the theorem that among the sub-Riemannian Manakov integrals \(g_{k,l}\) introduced in \eqref{sr_manakov_integrals}, there exist
			\[
			\frac{1}{2}\left(\dim\mathfrak{so}(n)-\rank\mathfrak{so}(n)\right)=\frac{1}{2}\left(\frac{n}{2}(n-1)-\left\lfloor\frac{n}{2}\right\rfloor\right)
			\]
			independent first integrals in involution, which defines an integrable system on the regular adjoint orbits of \(\mathfrak{so}(n)\). This completes the proof of Theorem~\ref{sr_thm}.\qed

		\begin{Remark}
			There is a straightforward generalisation which replaces the distribution of tangent vectors in \eqref{rolling_dist} with block matrices
			\[
			\Omega=\begin{pmatrix}
				0 & -\omega^T\\\omega & 0
			\end{pmatrix},
			\]
			where \(\omega\) is a \(k\times(n-k)\)-matrix, with sub-Riemannian length
			\[
			\sqrt{\Tr(\omega^T\hat{A}\omega)},
			\]
			for \(\hat{A}=\text{diag}(A_{n-k},\dots, A_n)\). When interpreted as a rolling distribution, this corresponds to the problem of rolling the (oriented) Grassmannian \(\widetilde{\mathbf{Gr}}(k,n)\) along its affine tangent space. This problem has been considered extensively before in \cite{leite,jurdjevic_conferencepaper} (see also the references therein) for the particular case of the standard Euclidean metric where \(\hat{A}\) may be taken to be the identity matrix. We remark however, that the result of Mishchenko and Fomenko used above does not work for \(k>1\). The reason for this is that, although Eq.~\eqref{important} holds for \(A=\textrm{diag}(0,\dots,0,A_{n-k},\dots,A_n)\) and \(B=\textrm{diag}(0,\dots,0,1,\dots,1)\), the matrix \(A\) is not regular, and so the proof in \cite{misch_fomenko} no longer applies.
			
		\end{Remark}
		
		\begin{Remark}
			It is somewhat frustrating that we must assume the \(A_2,\dots, A_n\) to be distinct in order to ensure that \(A\) remains regular. This frustration is shared with the rigid body, where the integrability results of Manakov, Mishchenko, and Fomenko both require the \(J_i\) to be distinct. In the event that these become equal, we obtain the so-called symmetric body. The proof of integrability breaks down as dependencies form between the Manakov integrals which were previously independent. On the other hand, the added symmetry of the body gives rise to additional conserved momenta. These, however, do not necessarily commute with one another, so one must consider the related notion of non-commutative integrability, as performed in \cite{singular_manakov,singular_manakov_2}. One can adapt their work to the sub-Riemannian case where some of the \(A_i\)s become equal. A related work is \cite{jsv} in which it is shown that specific sub-Riemannian metrics, which include our example for \(A_2=\dots=A_n\), are integrable in the non-commutative sense.
				
				Alternatively, there is a result due to Vinberg \cite{vinberg} which states that the limit of the algebra of integrals defined in Theorem~\ref{mf_thm} remains a commutative algebra with the same transcendence degree, even if \(A\) tends to a non-regular element. This can be used to establish complete integrability of the system as various \(A_i\)s converge to each other. We remark that the resulting integrable system is larger than the algebra of integrals defined in Theorem~\ref{mf_thm} and depends on the way in which the regular \(A\) converges to a non-regular element. These integrals and their dependence on the limit can be found explicitly by applying the work of Shuvalov~\cite{shuvalov}. 
		\end{Remark}

		\section{Concluding Remarks}

		A curious feature which distinguishes sub-Riemannian geometry from ordinary Riemannian geometry is the existence of singular geodesics. These are length minimising curves which, unlike normal geodesics, are not the projections of integral curves in the cotangent bundle for the Hamiltonian flow. In other words, these are geodesics that do not arise from solutions of the equations of motion, and are typically more difficult to find. Therefore, it is natural to ask what, if any, singular geodesics are present for the sub-Riemannian metric on \(\mathbf{SO}(n)\)? We note that for \(n=3\) the distribution is contact and so all geodesics are normal \cite[Proposition 4.38]{agrachev}.
		
		The relative equilibria of the standard Euler equations on \(\mathfrak{so}(3)\) correspond to the well-known steady motions of a rigid body around its principal axes. In our example for \(n=3\), the analogous solutions can be seen from Figure~\ref{rolling_balls} to correspond to motions where the ball rolls in a straight line across the table. For \(n>3\) a classification of all such relative equilibria is more challenging, and corresponds to finding critical points of the sub-Riemannian Hamiltonian on the adjoint orbits. This task has been considered for the case of geodesic flows in \cite{ratiu_equilibria}.

		Finally, we conclude with a conjecture. The famous Mishchenko-Fomenko conjecture (now a theorem \cite{sadetov}, see also \cite{bolsinov_MF}) asked whether every Lie group admits an integrable Riemannian metric. This was initially proven by Mishchenko~\&~Fomenko for semisimple Lie groups in \cite{misch_fomenko}. Our work motivates a related sub-Riemannian extension of this question: for any (semisimple) Lie group \(G\), does their exist an integrable sub-Riemannian metric of any given rank \(k\)? The known Riemannian case corresponds to \(k=\dim G\). In this paper we have showed that for \(G=\mathbf{SO}(n)\), the result is true for the distribution \(\mathcal{D}\) with rank \(k=n-1\).  The general 
		case would shed light on the class of systems discussed in \cite{bloch1994sub} and \cite{krakowski2015rolling}.

		\subsection*{Acknowledgements} 
		
		We wish to express our gratitude to the mathematics department at the University of Michigan for enabling this collaboration. We also thank Richard Montgomery for his helpful and constructive correspondence throughout the preparation of this paper and the referees for their useful reviews which greatly helped our exposition. This research was supported in part by NSF grant 2103026 and AFOSR grants FA9550-23-1-0215 and FA9550-23-1-0400

		\bibliographystyle{plain}
		\bibliography{sr_refs}

	\end{document}